\documentclass[12pt]{article}
\usepackage{amsmath}
\usepackage{amssymb}
\usepackage{vatola}

\def\q{\quad}
\def\qq{\qquad}
\def\mod#1{\ (\text{\rm mod}\ #1)}
\def\t{\hbox}
\def\({\left(}
\def\){\right)}
\def\qtq#1{\q\t{#1}\q}
\def\f{\frac}
\def\e{\equiv}

\def\b{\binom}

\let \pro=\proclaim
\let \endpro=\endproclaim
\begin{document}
\leftline{The paper will appear in International Journal of Number
Theory. } \par\q\par\q
 \centerline {\Large\bf
An extension of Stern's congruence}
$$\q$$

 \centerline{\bf \qq
Zhi-Hong Sun$^{1}$ and Lin-Lin Wang$^{2}$} \centerline{$^1$School of
Mathematical Sciences, Huaiyin Normal University,} \centerline{
Huaian, Jiangsu 223001, P.R. China} \centerline{E-mail:
zhihongsun@yahoo.com} \centerline{$^2$Center for Combinatorics,
 Nankai University,}
\centerline{Tianjin 300071, P.R. China,}  \centerline{E-mail:
wanglinlin$\underline{\q }$1986@yahoo.cn}

\abstract{Let $\{E_n\}$ be the Euler numbers. In the paper we
determine $E_{2^mk+b}-E_b$ modulo $2^{m+7}$, where $k$ and $m$ are
positive integers  and $b\in\{0,2,4,\ldots\}$.
\par\q
}
 \endabstract

\section*{1. Introduction}\par
Let $\Bbb N$ be the set of positive integers.
The Euler numbers
$\{E_n\}$ are given by
$$E_0=1,\q
E_{2n-1}=0,\q\sum_{r=0}^n\b{2n}{2r}E_{2r}=0\q(n=1,2,3,\ldots).$$ The
first few Euler numbers are shown below:
$$\align&E_0=1,\ E_2=-1,\ E_4=5,\ E_6=-61,\ E_8=1385,\ E_{10}=-50521,
\\&E_{12}=2702765,\ E_{14}=-199360981,\ E_{16}=19391512145.
\endalign$$
\par
For $k,m\in\Bbb N$ and $b\in\{0,2,4,\ldots\}$, in 1875 Stern ([2])
proved the following congruence, which is now known as Stern's
congruence:
$$E_{2^mk+b}\e E_b+2^mk\mod{2^{m+1}}.\tag 1.1$$
There are many modern proofs of (1.1). See for example [1,3,5,6].
\par
Let $b\in\{0,2,4,\ldots\}$ and $k,m\in\Bbb N$. In [4] the first
author showed that
$$\align &E_{2^mk+b}\e E_b+2^mk\mod{2^{m+2}}\qtq{for}m\ge 2,\tag 1.2\\& E_{2^mk+b}\e E_b+5\cdot
2^mk\mod{2^{m+3}}\qtq{for}m\ge 3\tag 1.3\endalign$$
 and
$$E_{2^mk+b}\e \cases
E_b+5\cdot 2^mk\mod{2^{m+4}}&\t{if $b\e 0,6\mod 8$,} \\ E_b-3\cdot
2^mk\mod{2^{m+4}}&\t{if $b\e 2,4\mod 8$}\endcases\qtq{for}m\ge
4.\tag 1.4$$
\par   From [4, Corollary 3.2] we know
that for any nonnegative integer $k$ and positive integer $n$,
$$\f{3^{2k+1}+1}4E_{2k}\e\sum_{r=0}^{n-1}(-1)^{n-1-r}
\b{k-1-r}{n-1-r}\b kr\f{3^{2r+1}+1}4E_{2r}\mod{2^{3n}}.\tag 1.5$$ In
the paper we use (1.5) to obtain a congruence for $E_{2^mk+b}-E_b$
modulo $2^{m+7}$, where $k,m\in\Bbb N$ and $b\in\{0,2,4,\ldots\}$.
In particular,
$$E_{2^mk+b}\e E_b+2^mk(7(b+1)^2-18)\mod{2^{m+7}}\qtq{for}m\ge 7.\tag 1.6$$
\par Throughout the paper we use $\Bbb Z_2$ to denote the set of rational
numbers whose denominator is odd.
\section*{2. The main result}
\pro{Lemma 2.1} Let $b\in\{0,2,4,\ldots\}.$ Then
$$\aligned E_b&
\e\cases 1-11b+15b^2+b^3-b^4\mod{2^{10}}&\t{if $4\mid b$},
\\289-91b-17b^2-7b^3+b^4\mod{2^{10}}&\t{if $4\mid
b-2$}.\endcases\endaligned$$\endpro Proof. By [4, Corollary 3.7],
$$\aligned E_b\e\cases\frac{1+172b-24b^2}{1-329b-74b^2-24b^3}\mod{2^{10}}
&\t{if $4\mid
b$},\\\frac{-7-308(b-2)+40(b-2)^2}{7+111(b-2)+102(b-2)^2-24(b-2)^3}\mod{2^{10}}
&\t{if $4\mid b-2$.}
\endcases\endaligned$$
It is easy to see that
$$\aligned &(1-329b-74b^2-24b^3)( 1-11b+15b^2+b^3-b^4)
\\&\ \e 1+172b-24b^2\mod{2^{10}} \endaligned$$ for $4\mid b$, and
that
$$\aligned&(7+111(b-2)+102(b-2)^2-24(b-2)^3)(289-91b-17b^2-7b^3+b^4)
\\&\ \e
-7-308(b-2)+40(b-2)^2\mod{2^{10}}\endaligned$$ for $4\mid b-2$. Thus
the result follows.

 \pro{Theorem 2.1} Let $b\in\{0,2,4,\ldots\}$ and $k,m\in\Bbb
N$. Then
$$\aligned&E_{2^mk+b}-E_b
\\&\e \cases 2k(-(b-7)^2+38+2k(3b-1+2k))\mod{2^{m+7}}
&\t{if $m=1$ and
$2\mid k$,}
\\2k(-(b+1)^2+6+2k(3b-1+2k))&
\\\qq-16(b+(-1)^{\f b2}) \mod {2^{m+7}}&\t{if $m=1$ and $2\nmid k$,}
\\
4k(7(b+1)^2-18&\\\qq +12k(b+1+4((-1)^{\f
b2}-k)))\mod{2^{m+7}}\q&\t{if $m=2$},
\\2^mk(7(b+1)^2-18
+2^mk(7-b))\mod {2^{m+7}}\q&\t{if $m\ge 3$.}
\endcases\endaligned$$

\endpro
 Proof. For $m=1,2,3,$ one can easily deduce the result from Lemma
 2.1. Now we assume that $m\geq 4.$  Set
$$f(k)=\f{3^{2k+b+1}+1}4E_{2k+b}\qtq{and}F(k)=2^{-3k}\sum_{r=0}^k\b
kr(-1)^rf(r).$$ From  [4, Lemma 2.3] we know that $F(k)\in\Bbb Z_2$.
By the binomial inversion formula we have $f(k)=\sum_{r=0}^k\b
kr(-2)^{3r}F(r)$ and so
$$\aligned
f(2^{m-1}k)&=F(0)+\sum_{r=1}^{2^{m-1}k}\b{2^{m-1}k}r(-2)^{3r}F(r)
\\&=F(0)+\sum_{r=1}^{2^{m-1}k}\f{2^{m-1}k}r\b{2^{m-1}k-1}{r-1}(-2)^{3r}F(r).\endaligned$$
For $r\ge 3$ we have $\f{2^{r-2}}r\in\Bbb Z_2$ and $F(r)\in\Bbb
Z_2$. Thus
$$\f{2^{m-1}\cdot 2^{3r}F(r)}r=2^{m+2r+1}\cdot\f{2^{r-2}}rF(r)\e
0\mod{2^{m+7}}.$$ Hence,
$$\aligned &\f{3^{2^mk+b+1}+1}4E_{2^mk+b}\\&=f(2^{m-1}k)\e
F(0)-2^{m+2}kF(1)+2^{m+4}k(2^{m-1}k-1)F(2)
\\&=\f{3^{b+1}+1}4E_b-2^{m-1}k(f(0)-f(1))-2^{m-2}k(f(0)-2f(1)+f(2))+2^{2m+3}k^2F(2)
\\&\e \f{3^{b+1}+1}4E_b-2^{m-2}k\(3f(0)-4f(1)+f(2)\)\mod{2^{m+7}}.\endaligned\tag 2.1$$
Putting $n=3$ and $k=b/2$ in (1.5) we see that
$$\aligned
\f{3^{b+1}+1}4E_b&\e\sum_{r=0}^2(-1)^{2-r}\b{b/2-1-r}{2-r}\b{b/2}r\f{3^{2r+1}+1}4E_{2r}
\\&=\f{(\f b2-1)(\f b2-2)}2-\f b2(\f b2-2)\cdot 7\cdot (-1)+
\f{\f b2(\f b2-1)}2\cdot 61\cdot 5
\\&=40b^2-84b+1\mod{2^9}.\endaligned$$
Thus,
$$\aligned &3f(0)-4f(1)+f(2)\\&\e
3(40b^2-84b+1)-4(40(b+2)^2-84(b+2)+1)+40(b+4)^2-84(b+4)+1
\\&=336-320b\e -176-64b\mod{2^9}.\endaligned\tag 2.2$$
 For $r\ge 3$ we see that $4^{r-3}/r\in\Bbb Z_2$ and so
 $$\b{2^{m-2}k}r{80}^r=2^{m+4+2r}\cdot\f{4^{r-3}}r\cdot k\cdot
 5^r\b{2^{m-2}k-1}{r-1}\e 0\mod{2^{m+9}}.$$
Thus, $$\aligned
3^{2^mk}&=(1+80)^{2^{m-2}k}=\sum_{r=0}^{2^{m-2}k}\b{2^{m-2}k}r{80}^r
\\&\e
1+5k\cdot 2^{m+2}+25k\cdot 2^{m+5}(2^{m-2}k-1)
\\&=1-195\cdot 2^{m+2}k+25k^2\cdot 2^{2m+3}\\&\e 1+61\cdot 2^{m+2}k
+2^{2m+3}k^2\mod{2^{m+9}}.\endaligned$$ This together with (2.1) and
(2.2) yields
$$\align&\f{3^{b+1}(1+61\cdot
2^{m+2}k+2^{2m+3}k^2)+1}4E_{2^mk+b}\\&\e\f{3^{2^mk+b+1}+1}4E_{2^mk+b}
\e\f{3^{b+1}+1}4E_b+2^{m-2}k(176+64b)\mod{2^{m+7}}.\endalign$$ Thus,
$$\f{3^{b+1}+1}4\(E_{2^mk+b}-E_b\)\e  2^mk(2^{m+1}k-55) 3^bE_{2^mk+b}+2^mk(44+16b)\mod{2^{m+7}}.\tag 2.3$$ By [3,
Corollary 7.9] or Lemma 2.1,
$$E_b\e\cases 3b^2-11b+1\mod{2^7}&\t{if}\q 4\mid
b,\\b^2-23b+41\mod{2^7}&\t{if}\q 4\mid b-2.\endcases\tag 2.4$$
Replacing $b$ with $2^mk+b$ in (2.4) we get
$$E_{2^mk+b}\e\cases 3b^2-11b+1+2^mk(6b-11)\mod{2^7}&\t{if}\q 4\mid
b,\\b^2-23b+41+2^mk(2b-23)\mod{2^7}&\t{if}\q 4\mid b-2.\endcases\tag
2.5$$ As
$$\aligned 3^b&=(1+8)^{b/2}=1+\b{b/2}18+\b{b/2}28^2+\cdots\e
8b^2-12b+1\\&\e\cases -12b+1\mod{2^7}&\t{if $4\mid b$,}
\\20b-31\mod{2^7}&\t{if $4\mid b-2$,}\endcases\endaligned$$
 by (2.5) we have
$$3^bE_{2^mk+b}\e\cases (3b^2-11b+1+2^mk(6b-11))(-12b+1)\mod{2^7}
&\t{if}\q 4\mid b,\\(b^2-23b+41+2^mk(2b-23))(20b-31) \mod{2^7}&
\t{if}\q 4\mid b-2.\endcases$$ If $4\mid b$, then
$$(3b^2-11b+1)(-12b+1)=-36b^3+135b^2-23b+1\e 7b^2-23b+1\mod{2^7}$$
and
$$(6b-11)(-12b+1)\e -11\e -3\mod 8.$$
If $4\mid b-2$, then $32\mid (b-2)(b+2)$ and so
$$\align (b^2-23b+41)(20b-31)
&=((b-2)^2-19(b-2)-1)(20(b+2)-71)
\\&\e -20(b+2)-71(b^2-23b+41)
\\&\e -7b^2+13b-7\mod{2^7}\endalign$$ and
$$(2b-23)(20b-31)=(2(b-2)-19)(20b-31)\e 19\cdot 31\e -3\mod 8.$$
Thus
$$3^bE_{2^mk+b}\e\cases 7b^2-23b+1-3\cdot 2^mk\mod{2^7}&\t{if}\q 4\mid
b,\\ -7b^2+13b-7-3\cdot 2^mk \mod{2^7}& \t{if}\q 4\mid
b-2.\endcases$$
 Substituting this into (2.3) we obtain
$$\aligned&\f{3^{b+1}+1}4(E_{2^mk+b}-E_b)\\&\e\cases 2^mk(-b^2+b-11-2^mk)
\mod{2^{m+7}}&\t{if}\q 4\mid b,\\2^mk(b^2+5b+45+3\cdot 2^mk)
\mod{2^{m+7}}&\t{if}\q 4\mid b-2.\endcases\endaligned\tag 2.6$$ It
is easily seen that
$$ \f{3^{b+1}+1}4=\f{3(1+8)^{b/2}+1}4\e\f{3(1+\f b2\cdot 8+\b{b/2}2\cdot
8^2)+1}4=6b^2-9b+1\mod{2^7}.$$
\par If $4\mid b$, then $b^3=4^3(\f b4)^3\e 4^3\cdot \f
b4=16b\mod{2^7}$. One can easily see that
$$(6b^2-9b+1)(3b^2-7b+1)\e 1\mod{2^7},\q 3b^2-7b+1\e b+1\mod 8$$
and
$$(-b^2+b-11)(3b^2-7b+1)\e -b^2+14b-11\mod{2^7}.$$
Thus, by (2.6) we have
$$\aligned E_{2^mk+b}-E_b&
\e 2^mk(-b^2+b-11-2^mk)(3b^2-7b+1)
\\&\e 2^mk(-b^2+14b-11-(b+1)\cdot 2^mk)
\\&\e 2^mk(7b^2+14b-11+2^mk(7-b))\mod{2^{m+7}}.
\endaligned$$
\par If $4\mid b-2$, then $(b+2)^3\e 16(b+2)\mod{2^7}$.
One can easily see that
$$\align &(6b^2-9b+1)(-3(b+2)^2-7(b+2)+3)\e 1\mod{2^7},
\\&-3(b+2)^2-7(b+2)+3\e b+5\mod 8\endalign$$
and
$$\align &(b^2+5b+45)(-3(b+2)^2-7(b+2)+3)\\&
=((b+2)^2+(b+2)+39)(-3(b+2)^2-7(b+2)+3) \\& \e -b^2+14b+21\e
7b^2+14b-11\mod{2^7}.\endalign$$ Thus, by (2.6) we have
$$\align E_{2^mk+b}-E_b&
\e 2^mk(b^2+5b+45-2^{m+1}k)(-3(b+2)^2-7(b+2)+3)
\\&\e 2^mk(7b^2+14b-11+ 2^mk(7-b))\mod{2^{m+7}}.
\endalign$$
 Note that $7b^2+14b-11=7(b+1)^2-18$. Combining all
the above we prove the theorem.
\par\q

 \pro{Corollary 2.1} Let
$k,m\in\Bbb N$ with $m\ge 2$. Then
$$E_{2^mk}\e\cases 4k(-48k^2+60k-11)+1\mod{2^{m+7}}&
\t{if $m=2$},\\ 2^mk(-11+7\cdot 2^mk)+1\mod{2^{m+7}}&\t{if $m\geq
3$}\endcases$$ and
$$E_{2^mk+2}\e\cases 4k(-48k^2-12k+45)-1\mod{2^{m+7}}&\t{if $m=2$},
\\2^mk(45+5\cdot 2^mk)-1\mod{2^{m+7}}&\t{if $m\geq 3$}.\endcases$$
\endpro

\pro{Corollary 2.2} Let $b\in\{0,2,4,\ldots\}$ and $k,m\in\Bbb N$
with $m\ge 3$. Then
$$E_{2^mk+b}\e
E_b-2^mk((b+1)^2+10+2^mk(b+1))\mod{2^{m+6}}.$$
\endpro
Proof. Since $7(b+1)^2-18\e-(b+1)^2-10\mod{2^6}$ and
$7-b\e-b-1\mod{8}$, the result follows from Theorem 2.1.
\pro{Corollary 2.3} Let $b\in\{0,2,4,\ldots\}$ and $k,m\in\Bbb N$
with $m\ge 5$. Then
$$E_{2^mk+b}\e\cases
E_b-11\cdot 2^mk\mod{2^{m+5}}&\t{if}\q b\e 0,14\mod{16}, \\
E_b+13\cdot 2^mk\mod{2^{m+5}}&\t{if}\q b\e 2,12\mod{16},
\\
E_b-3\cdot 2^mk\mod{2^{m+5}}&\t{if}\q b\e 4,10\mod{16},
\\
E_b+5\cdot 2^mk\mod{2^{m+5}}&\t{if}\q b\e 6,8\mod{16}
.\endcases$$\endpro Proof. By Corollary 2.2 we have
$$E_{2^mk+b}\e
E_b-2^mk((b+1)^2+10)\mod{2^{m+5}}.$$ Thus the result follows.
\par In conclusion we pose the following conjecture.
\pro{Conjecture 2.1} Let $m,n\in\Bbb N$, $m\ge n$ and
$b\in\{0,2,4,\ldots, 2^{m+n-1}-2\}$. Then
$$E_{2^mk+b}-E_b\e
E_{2^mk+2^{m+n-1}-2-b}-E_{2^{m+n-1}-2-b}\mod{2^{m+n}}.$$
\endpro

\end{document}